\documentclass[12pt]{amsart}

\usepackage{latexsym}
\usepackage[all]{xy}
\usepackage{amsmath, amsthm }
\usepackage{amssymb}
\usepackage{amsfonts}
\usepackage{color}
\usepackage{mathtools}
\usepackage[T1]{fontenc}

\usepackage{hyperref}
\usepackage{xcolor}
\hypersetup{colorlinks=true,urlcolor=blue,linkcolor=blue,citecolor=blue}

\usepackage{cite}

\newcommand {\QCoh} {\mathsf{QCoh}}

\newcommand{\csMod}{\mathsf{csMod}}
\newcommand{\sCR}{\mathsf{sCR}}
\newcommand{\Cinf}{\widehat{C}_{\mathrm{inf}}}
\newcommand{\CDR}{\widehat{C}_{DR}}
\newcommand{\cL}{\mathcal{L}}

\newcommand {\Spf} {\mathsf{Spf}}

\newcommand {\uMap} {\underline{\mathbf{Map}}}

\newcommand {\OO} {\mathcal{O}}

\newcommand {\Fol}{\mathcal{F}ol}

\newcommand {\F} {\mathcal{F}}

\newcommand {\Ga} {\mathbb{G}_a}
\newcommand {\Gm} {\mathbb{G}_m}
\newcommand {\T} {\mathbb{T}}

\newcommand{\QQ}{\mathbb{Q}}
\newcommand{\RR}{\mathbb{R}}
\newcommand{\VV}{\mathbb{V}}

\newcommand{\ZZ}{\mathbb{Z}}
\newcommand{\LL}{\mathbb{L}}
\newcommand{\cH}{\mathcal{H}}
\newcommand {\Spec} {\mathbf{Spec}}
\newcommand  {\dg}     {\mathbf{dg}}
\newcommand  {\csCR}     {\mathbf{csCR}}
\newcommand  {\fdg}     {\mathbf{dg}^{fil}}

\newcommand  {\mdg}     {\mathbf{dg}^{gr}}

\newcommand  {\medg}     {\epsilon-\mathbf{dg}^{gr}}
\newcommand {\scat} {\mathbf{Cat}_{\infty}}

\newcommand  {\dAff}     {\mathbf{dAff}}
\newcommand  {\dChAff}     {\mathbf{dChAff}}

\newcommand  {\dSt}   {\mathbf{dSt}}

\newcommand{\s}{\infty}

\newcommand{\D}{\mathcal{D}}

\newcommand{\uHom}{\underline{Hom}}

\theoremstyle{plain}
\newtheorem{thm}{Theorem}[section]
\newtheorem{df}[thm]{Definition}
\newtheorem{prop}[thm]{Propositon}
\newtheorem{rmk}[thm]{Remark}
\newtheorem{cor}[thm]{Corollary}

\newtheorem{lem}[thm]{Lemma}

\newtheorem{quest}[thm]{Question}

\begin{document}

\title{Infinitesimal derived foliations}
\author{Bertrand To\"en and Gabriele Vezzosi}

\maketitle

\begin{abstract} 
We introduce a notion of \emph{infinitesimal derived foliation}. We prove it is
related to the classical notion of infinitesimal cohomology, and satisfies some formal 
integrability properties. We also provide some hints on how infinitesimal derived foliations 
compare to 
our previous notion of derived foliations of \cite{derfol,derfolbook}.
\end{abstract}

\tableofcontents

\section*{Introduction}

The object of this short note is to present a notion of \emph{infinitesimal derived foliations}
on general schemes and stacks over bases of arbitrary characteristics (see 
Definition \ref{dinfol}). This notion 
is slightly different from the one exposed in \cite{derfol} but it is based on similar ideas. By definition, 
an infinitesimal derived foliation on $X$ consists of a perfect complex $E$ (the cotangent
complex of the foliation) together with an action of $G_0=\Omega_0\Ga$, the loop group
of the additive group, on the total space $\VV(E[-1])$ of the shift $E[-1]$, compatible with the natural gradings.
This is very similar to the definition in \cite{derfol}, except that the group $BKer$ is replaced
here by $G_0$, and the shift by $1$ is replaced by a shift by $-1$. These changes may sound
purely aesthetic choices, but outside of characteristic zero these two notions of derived
foliations are, in fact, different.

A major advantage of infinitesimal derived foliations is that they satisfy several 
formal integrability results (see Corollary 
\ref{cint}) that do \emph{not} hold for derived foliations. The key point here is that 
the natural cohomology theory related to infinitesimal derived foliations is 
(Hodge-completed derived) infinitesimal 
cohomology, whereas derived foliations are related to (Hodge-completed derived) 
de Rham cohomology. This
relation with infinitesimal cohomology makes it possible to formally integrate any 
(smooth) infinitesimal derived foliation by a smooth formal groupoid. In a way, 
infinitesimal derived foliations can be thought of as derived foliations endowed with extra 
conditions/structures. This resembles the difference between foliations 
and \emph{restricted foliations} in characteristic $p>0$ (see for instance \cite[Prop. 2.3]{MR927978} as well as \cite[Def. 3.6]{MR927960}).

The first section of this paper contains results on \emph{derived affine stacks}, that are a derived extension of affine stacks
of \cite{chaff}. We present in particular algebraic models for those, using a certain model 
category of cosimplicial-simplicial commutative algebras, which might be of independent interests.
These models will be particularly useful in our last section where we construct the red-shift
construction for graded derived affine stacks (see \S 5). The second section presents the general
definition of infinitesimal derived foliations and their basic properties. In section $\S 3$ we study
the cohomology of infinitesimal derived foliations and prove that it recovers
the usual notion of infinitesimal cohomology when applied to the tautological foliation. 
This cohomology comparison is then used to construct, for any 
smooth infinitesimal derived foliation $\F$ on a smooth scheme $X$, a smooth formal groupoid $G_\F \to X$, 
and we show that $\F \mapsto G_\F$ produces an equivalence of categories. This result has to be
related with \cite[Prop. 2.3]{MR927978}, and in a way, is one possible generalization of Ekedahl's result. Finally, in the
very last section, we study the relations between infinitesimal derived foliations and our original 
notion of derived foliations (see \cite{derfol,derfolbook}). Unfortunately the results here are 
not optimal, and we leave several open questions on the way. We think that there should exists
a forgetful functor from infinitesimal derived foliations to derived foliations, obtained by 
some application of the red-shift construction, but we could not find a completely satisfying construction. The reader
should consider this last part as a suggestion for future research. 

As a final comment, this note will appear in an extended form as a chapter of a book in 
preparation on the subject (\cite{derfolbook}). \\

\textbf{Acknowledgments.} We thank J. Nuiten and J. Fu for several discussions on the
content of this paper. The idea of infinitesimal derived foliations has been inspired by 
J. Fu's current work on the structure of the Chevalley-Eilenberg complex of partition Lie algebras.

\section{Derived affine stacks}

\subsection{Cosimplicial-simplicial algebras}

We let denote by $sCR$ be the category of simplicial commutative rings. It is endowed with
its usual \emph{model category} structure for which equivalences and fibrations are defined on the
underlying simplicial sets.

Similarly, we denote by $csCR$ the category of cosimplicial objects in $sCR$, or equivalently 
of cosimplicial-simplicial commutative rings. Its objects will be denoted by $A^*_*$, where the
superscript will refer to the cosimplicial direction and subscript to the simplicial direction.
The category $csCR$ can be endowed with a \emph{tensored} and \emph{cotensored} structure over the category $sSet$
of simplicial sets. We warn the reader here that there are several natural such structures and that 
the choice we make here is important for the sequel. 
For a simplicial set $K$ and an object $A$
of $csCR$, we define the cotensored structure by 
$A^K$ by the formula
$$(A^K)^q_p:=(A^q_p)^{K_q}$$
for $[p] \in \Delta^{op}$ and $[q] \in \Delta$. 
The tensored structure $K \otimes A$, for $A \in csCR$ and $K \in sSets$, is defined by the usual adjunction formula
$Hom(K\otimes A,B) \simeq Hom(A,B^K)$, for arbitrary $B \in csCR$. Similarly, we have a canonical simplicial enrichment
for which the simplicial Homs are denoted by $\uHom$ and are defined by 
$$\uHom(A,B)_p:=Hom_{csCR}(A,B^{\Delta^p}).$$

By considering the natural embedding $Mod_\ZZ \subset C(\ZZ)$, from $\ZZ$-modules to cochain complexes, 
any object $A \in csCR$ can be considered, by forgetting the multiplicative structure, 
as a cosimplicial-simplicial object in $C(\ZZ)$. We denote by $Tot^{\pi}(A)$
the complex equals to $\prod_{q-p=n}A^q_p$ in cohomological degree $n$
and whose differential is given by the alternated sums of the faces and cofaces of the
cosimplicial-simplicial diagram. Equivalently, $A$ can be first turned into a 
bi-complex, by individual totalisations along the cosimplicial and the simplicial directions, and
$Tot^\pi(A)$ can then be realized as the product-total complex of this bicomplex. 
Note also that $Tot^\pi(A)$ is a model for the homotopy limit-colimit of $A$
$$Tot^\pi(A) \simeq \mathrm{holim}_{[p] \in \Delta}\left(\mathrm{hocolim}_{[q] \in \Delta^{op}} A^q_p \right) \in C(\ZZ).$$

\begin{df}
A morphism $A \to B$ in $csCR$ is a \emph{completed quasi-isomorphism}
if the induced morphism $Tot^\pi(A) \to Tot^\pi(B)$ is a quasi-isomorphism of complexes.
\end{df}

We warn the reader that even though the categories $csCR$ and $scCR$, of cosimplicial-simplicial
commutative rings and of simplicial-cosimplicial rings, are equivalent, the notion 
of completed quasi-isomorphisms is not the same as the notion of weak equivalences used
for instance in \cite{bracanu}.

\begin{thm}\label{tsome}
The category $csCR$ is endowed with a model category structure for which the weak equivalences
are the completed quasi-isomorphism, and the fibrations are the epimorphisms. This model 
category structure is moreover cofibrantly generated and is a simplicial model category for
the above mentioned simplicial enrichment.
\end{thm}

\textit{Proof.} As fibrations\footnote{Note that $f:A\to B$ in $csCR$ is surjective iff $Tot^\pi(f)$ is.} and equivalences are defined via the functor
$Tot^\pi$, the proof of the theorem consists of transferring the projective model structure 
on $C(\ZZ)$ to $csCR$ via the functor $Tot^\pi$. For this we apply the path object
argument (see for instance \cite[\S 2.6]{MR2016697}).

The functor $Tot^\pi$ is a right adjoint whose left adjoint sends a complex $E \in C(\ZZ)$ to 
the free commutative ring over the cosimplicial-simplicial module given by
$$(p,q) \mapsto (C^*(\Delta^p) \otimes C_*(\Delta^q) \otimes E)^0/Im\, d^{-1},$$
where $C^*(\Delta^p)$ and $C_*(\Delta^q)$ are the cohomology and homology complexes
of the standard simplex $\Delta^p$ and $\Delta^q$. This left adjoint 
will be denoted by $L\phi$, as being the composition of 
the free commutative ring functor $L : csMod \to csCR$, with the
functor $\phi : C(\ZZ) \to csMod$ defined by the above formula.
Clearly, the functor $Tot^\pi$ preserves small objects (but does not 
commute with filtered colimits in general). In order to apply \cite[\S 2.6]{MR2016697} we 
thus simply have to prove that $csCR$ possesses a fibrant replacement functor
and that all fibrant object possesses a path object. But all objects
are fibrant by definition, and we are thus reduced to show the existence
of path objects. For this we use the following lemma.

\begin{lem}\label{l1}
For any simplicial set $K$ and any $A \in csCR$, there exists a canonical 
quasi-isomorphism
$$Tot^\pi(A^K) \simeq Tot^\pi(A)^K,$$
where for $E \in C(\ZZ)$ and $K \in sSet$ we define
$$E^K :=\mathrm{holim}_{[p] \in \Delta} E^{K_p}.$$
\end{lem}

\textit{Proof of the lemma.} We remind that for any functor $F : \Delta \times \Delta \to \dg$
(or more generally any bi-cosimplicial object in a complete $\s$-category), 
the homotopy limit of $F$ can be computed, up to a quasi-isomorphism, 
either component-wise, or diagonally
$$\mathrm{holim}_{[p] \in \Delta} F(p,p) \simeq \mathrm{holim}_{[p]\in \Delta} \mathrm{holim}_{[q] \in \Delta} F(p,q).$$
We apply this to the functor sending $([p],[q])$ to the complex $N(A^q_*)^{K_p}$, where 
$N : sMod \to dg$ is the normalization construction. Computing the homotopy limit
diagonally yields $Tot^\pi(A^K)$, whereas the double homotopy limit construction yields $Tot^\pi(A)^K$.
\hfill $\Box$ \\

Lemma \ref{l1} now implies the existence of path objects. Indeed, for $A \in csCR$ an explicit
path object is given by $A \to A^{\Delta^1} \to A^{\partial \Delta^1}= A \times A$. The fact that
the constant map morphism $A \to A^{\Delta^1}$ is a completed quasi-isomorphism follows from
the lemma \ref{l1}. 
The fact that $A^{\Delta^1} \to A^{\partial \Delta^1}= A \times A$ is a fibration simply follows from the fact that
it is a levelwise epimorphism as this can be checked easily using the explicit formula for 
$A^K$. 

In order to finish the proof of the theorem, we are left to showing that the simplicial enrichement
is compatible with the model category structure. In other words, we must show that for
any cofibration of simplicial sets $i : K \hookrightarrow L$ and any fibration $f : A \to B$ in $csCR$, 
the induced morphism
$$f^i : A^L \longrightarrow B^L \times_{B^K}A^K$$
is a fibration, which is also an equivalence if $i$ is a trivial cofibration or $f$ is a trivial
fibration. The fact that the morphism $f^i$ is a fibration is clear because fibrations are epimorphisms
and because of the explicit formula for the exponentiation by a simplicial set in $csCR$. Finally, 
the fact that $f^i$ is also a weak equivalence when $i$ or $f$ is so simply follows from 
Lemma \ref{l1}.
\hfill $\Box$ \\

\begin{df}
The \emph{$\s$-category of cs-rings} is the $\s$-category 
obtained from $csCR$ by inverting the completed quasi-isomorphisms. It is denoted by
$\csCR$ (or by $\csCR_\mathbb{U}$ if one wants to specify cosimplicial-simplicial rings belonging
to a given universe $\mathbb{U}$).
\end{df}

\begin{rmk}
\emph{We note also that the model structure on $csCR$ can be equivalently be
obtained by first constructing a model structure on $csMod_\ZZ$, the category of cosimplicial-simplicial 
abelian groups, where again the equivalences and fibrations are defined via
the $Tot^\pi$ functor. The model structure on $csCR$ can then be obtained by
transferring along the forgetful functor $csCR \to csMod_\ZZ$, with left adjoint
given by the free commutative cosimplicial-simplicial ring.}
\end{rmk}

\begin{rmk}\emph{To finish, we mention that it is possible to a construct an $\s$-functor
$Tot^\pi : \csCR \longrightarrow LSym-Alg$, from our $\s$-category of cosimplicial-simplicial
commutative rings to the $\s$-category of \emph{derived rings}, namely modules over the
$LSym$-monad. Here we denote by $LSym$ the monad
associated to the derived commutative operad as in \cite[Ex. 3.71]{bracanu}. This $\s$-functor simply sends $A$ to $\lim_{q}A^q_*$, where
the limit is taken in $LSym-Alg$ and $A^q_*$ is considered as a connective $LSym$-algebra and thus
as an object in $LSym-Alg$. This $\s$-functor is the right adjoint of an adjunction, and
preserves free objects over perfect complexes. It is likely that $Tot^\pi$ is an equivalence of $\s$-categories
when restricted to nice enough objects.}
\end{rmk}

\subsection{Spectrum of cosimplicial-simplicial rings}

For a fixed commutative simpicial ring $k \in sCR$, we can 
work relatively over $k$ and define the model category $k-csCR$ of cosimplicial-simplicial 
commutative $k$-algebras, and its associated $\s$-category $k-\csCR$. As usual, we have a 
natural equivalence of $\s$-categories
$$k-\csCR \simeq k/\csCR,$$
where $k$ is considered as an object in $\csCR$ which is constant is the cosimplicial direction.

We consider $dAff_k:=(k-sCR)^{op}$, 
the model category of derived affine $k$-schemes, which is defined to be the opposite category of
that of simplicial commutative $k$-algebras. Remind from \cite{hagII} that it can be endowed
with the fpqc model topology, and that we can consider the model category of (hyper-complete) stacks
$dAff_k^{\sim,fpqc}$. There are here some set-theoretical issues, that can be solved,
as usual, by fixing two universes $\mathbb{U} \in \mathbb{V}$. By definition $k-sCR$ refers here
to $\mathbb{U}$-small simplicial $k$-algebras, and $dAff_k^{\sim,fpqc}$ is then the category of
functors $k-sCR \to sSet_\mathbb{V}$ to $\mathbb{V}$-small simplicial sets.

We then consider a functor
$$\Spec^\Delta : k-csCR_\mathbb{V}^{op} \longrightarrow dAff_k^{\sim,fpqc},$$
from $\mathbb{V}$-small cosimplicial-simplicial commutative $k$-algebras to $dAff_k^{\sim,fpqc}$. 
It is defined by 
sending $A \in k-csCR_\mathbb{V}$ to the functor $\Spec^\Delta\, A : sCR \to sSet_\mathbb{V}$ defined by 
$$\Spec^\Delta\, A : B \mapsto \uHom(A,B)= Hom_{csCR}(A,B^{\Delta^\bullet}),$$
where $\uHom$ are the simplicial $Hom$'s of the simplicial enrichement described in the 
previous section, and $B$ is considered as an object in $k-csCR_\mathbb{V}$ by viewing it as constant
in the cosimplicial direction.

\begin{prop}\label{p1}
The functor $\Spec^\Delta : k-csCR_\mathbb{V}^{op} \longrightarrow dAff_k^{\sim,fpqc}$
is right Quillen, and the induced functor on $\s$-categories
$$\Spec^\Delta : k-\csCR_\mathbb{U} \longrightarrow \dSt_k$$
is fully faithful. The essential image of $\Spec^\Delta$ is the smallest
full sub-$\s$-category of $\dSt_k$ containing the objects
$K(\Ga,n)$ for various $n$, and which is stable by
 $\mathbb{U}$-small limits.
\end{prop}

\textit{Proof.} This is proven in a very similar manner than \cite[Cor. 2.2.3]{chaff}.
The left adjoint to $\Spec^\Delta$, denoted by $\OO$, sends 
a representable presheaf $h^B=Hom(B,-) : k-sCR \to sSet$ to 
$B \in k-csCR$, which is considered as constant in the cosimplicial direction. It also sends
an object of the form $K \times h^B$, for a simplicial set $K \in sSet$, to 
$B^K \in k-csCR$. Finally, it is uniquely defined by these properties together with the requirement that 
it sends colimits in $dAff_k^{\sim,fpqc}$ to limits in $k-csCR$. 

To prove that $\Spec^\Delta$ is right Quillen we use that $dAff_k^{\sim,fpqc}$ is a left Bousfield localization
of the levelwise projective model structure on the category of simplicial presheaves 
$Fun(k-sCR_\mathbb{U},sSet)$ obtained by inverting fpqc-local equivalences and 
equivalences in $sCR$ (see \cite{hagII}). Therefore, by general facts about Bousfield localizations,
it is enough to show that $\Spec^\Delta$ is a right Quillen functor for the
levelwise projective model structure, and moreover that for any cofibrant $A \in k-csCR$, 
$\Spec^\Delta\, (A)$ is a fibrant object in $dAff_k^{\sim,fpqc}$. The first of these statements 
easily follows from the fact that $k-csCR$ is a simplicial model category in which 
every object is fibrant, and thus
that if $A \to A'$ is a (trivial) cofibration in $k-csCR$, for any $C \in k-sCR$ the morphism
$\uHom(A',C) \to \uHom(A,C)$ 
is a (trivial) fibration of simplicial sets. For the second of these statements, let $A \in k-csCR$ be 
a cofibrant object. Any equivalence $B \to B'$ in $k-sCR$ obviously induces an equivalence
in $k-csCR$ when considered both $B$ and $B'$ as constant in the cosimplicial direction, and thus
$\uHom(A,B) \to \uHom(A,B')$ is an equivalence of simplicial sets. This shows that 
$\Spec^\Delta(A)$ is fibrant when the topology is the trivial topology. Finally, 
if $B=\mathrm{holim}_i B_i$ is a homotopy limit in $k-sCR$, it is also a homotopy limit in 
$k-csCR$, and thus the natural morphism $\uHom(A,B) \to \mathrm{holim}_i \uHom(A,B_i)$
is an equivalence. As the fpqc model topology is subcanonical, this implies that 
$\Spec^\Delta(A)$ is indeed a stack for the fpqc topology, and thus fibrant 
as an object in $dAff_k^{\sim,fpqc}$ (see \cite{hagII}).

To finish the proof of the proposition, let $A \in k-csCR_\mathbb{U}$ and 
$X=\Spec^\Delta(A)$. By definition of the simplicial structure on $k-csCR$, the functor
$X$ is the geometric realization of the simplicial object $q \mapsto h^{A_*^q}$, where
$h^{A_*^q} : k-sCR \to sSet$ is corepresented by $A_*^q \in k-sCR$. In other words, we have 
$\Spec^\Delta(A) \simeq \mathrm{hocolim}_{q} h^{A^q_*}$. Therefore, the adjunction morphism
$$A \longrightarrow \OO(\Spec^\Delta(A))$$
is the canonical morphism $A \to \mathrm{holim}_{q}A^q_*$ in $Ho(k-csCR_\mathbb{U})$. This canonical 
morphism is obviously an equivalence, showing the fully faithfulness property in the proposition.
It is then formal that the essential image of $\RR\Spec^\Delta$ is stable by limits. 
It also contains the objects $K(\Ga,n)$, as these are the images of $L\phi(\ZZ[-n])$, where
$L\phi : C(\ZZ) \to k-csCR$ is the left adjoint to $Tot^\pi$. Finally, 
as $k-csCR$ is cofibrantly generated, any object is equivalent to a $I$-cell object, where
$I$ is the image by $L\phi$ of the generating cofibrations in $C(\ZZ)$. The images of 
$I$ by $\Spec^\Delta$ are equivalent to morphisms of the form $* \to K(\Ga,n)$, so that 
 any object is the essential image of $\RR\Spec^\Delta$ lies in the smallest
sub-$\s$-category containing the $K(\Ga,n)$ and stable by limits.
\hfill $\Box$ \\

\begin{df}
The $\s$-category of \emph{($k$-linear) derived affine stacks} is the essential image of 
$\RR\Spec^\Delta : k-\csCR_\mathbb{U}^{op} \hookrightarrow \dSt_k$. It is denoted by 
$\dChAff_k$.
\end{df}

The above notions and results also have \emph{graded versions}, for which the details are left to the reader.
We denote by $k-csCR^{gr}$ the category of graded cosimplicial-simplicial commutative $k$-algebras. We endow
this category with the model category structure for which equivalences and fibrations
are defined by the graded Tot functor to graded complexes $Tot^\pi : k-scCR^{gr} \to C(\ZZ)^{gr}$, 
defined by taking the $Tot^\pi$ of each graded component. The proof of the existence of this model 
category structure follows the same lines as in the non-graded case. The corresponding 
$\s$-category will be denoted by $k-\csCR^{gr}$.

We consider the strict multiplicative group $\Gm^{st} \in \dAff_k^{\sim,fpqc}$, given by the
functor sending $B \in k-sCR$ to the group $B_0^*$ of invertible elements in the ring $B_0$
of $0$-simplices in $B$. This is a group object in the model category $\dAff_k^{\sim,fpqc}$, which 
is not a fibrant object (it does not preservers equivalences in $B$), but its image
in $\dSt_k$ is the usual multiplicative group scheme $\Gm$. Indeed, 
$\Gm^{st}=h^{k[t,t^{-1}]}$ is corepresented by $k[t,t^{-1}]$, and,
as shown in \cite{hagII}, a fibrant model for $h^A$ is the derived scheme $\Spec\, A : B \mapsto Map(A,B)$.
In particular, we have that the model category 
$\Gm^{st}-\dAff_k^{\sim,fpqc}$, of $\Gm^{st}$-equivariant objects in $\dAff_k^{\sim,fpqc}$, is a model
for the $\s$-category $\Gm-\dSt_k$, of $\Gm$-equivariant derived stacks, called
\emph{graded derived $k$-stacks}.

The graded version of the spec functor is the functor
$$\Spec^{\Delta,gr} : (k-csCR^{gr})^{op} \longrightarrow \Gm^{st}-\dAff_k^{\sim,fpqc},$$
sending $A \in k-csCR^{gr}$ to $\uHom(A,-)$, endowed with the natural
$\Gm^{st}$-action coming from the grading on $A$. More explicitly, 
for $B \in k-sCR$, the set of q-simplices of $\Spec^{\Delta,gr}(A)(B)$ is the set
of morphisms $Hom(A^q,B)$, which is endowed with a $B_0^*$-action as follows. 
For $b \in B_0^*$ and $f : A^q \to B$, we define $bf : A^q \to B$ by the formula
$(bf)(x)=b^n.f(x)$ for a homogeneous element $x$ of degree $n$. As in the non-graded
case, $\Spec^{\Delta,gr}$ is a right Quillen functor, and the induced $\s$-functor 
$$\RR\Spec^{\Delta,gr} : (k-\csCR^{gr}_\mathbb{U})^{op} \longrightarrow \Gm-\dSt_k$$
is fully faithful. Its essential image is the smallest sub-$\s$-category 
containing all the objects $K(\Ga^{\chi},n)$, for $n\geq 0$ and $\chi \in \ZZ$,  
where $\Gm$ acts on $\Ga$ with weight $\chi \in \ZZ$, 
and which is stable by $\mathbb{U}$-small limits. 

\begin{df}
The $\s$-category of \emph{graded derived affine $k$-stacks} is the essential image of 
$\RR\Spec^{\Delta,gr} : (k-\csCR_\mathbb{U}^{gr})^{op} \hookrightarrow \Gm-\dSt_k$. It is denoted by 
$\dChAff_k^{gr}$.
\end{df}

We finish this section by recalling the notion of \emph{linear derived stack}, already introduced and studied in \cite{mon}. 
We first notice that the Tot functor
$Tot^\pi : csMod_k \longrightarrow C(\ZZ)$ admits a natural lax monoidal structure. Therefore, 
it gives rise to an $\s$-functor
$$Tot^\pi : \csMod_k \longrightarrow \dg_k,$$
where $\csMod_k$ is the $\s$-category of cosimplicial-simplicial $k$-modules.
This $\s$-functor turns out to be an equivalence of $\s$-categories, and we denote by
$\phi : \dg_k \to \csMod_k$ its inverse.

For any object $E \in \dg_k,$ we consider
$\phi(E) \in \csMod_k$, and define
$$Sym^{\Delta}(E):=L\phi(E) \in csCR^{gr}$$
the free cosimpicial-simplicial commutative ring generated by $\phi(E)$ (same notation as in the proof of Theorem \ref{tsome}). Being a free
commutative ring, $Sym^{\Delta}(E)$ comes natural equipped with a graduation, and thus
is considered as an object in the $\s$-category $\csCR^{gr}$. 

\begin{df}
The \emph{$k$-linear derived stack associated to a complex $E \in \dg_k$}, 
is the graded derived affine $k$-stack defined by
$$\VV(E):=\mathbb{R}\Spec^{\Delta,gr}\,(Sym^{\Delta}(E)) \in \Gm-\dSt_k.$$
\end{df}

Note that by construction the functor of points of $\VV(E)$ sends $B \in k-sCR$
to $Map_{\dg_k}(E,N(B))$, where $N(B)$ is the normalisation of $B$ as a simplicial
$k$-module. The $\Gm$-action is then the natural action of invertible elements 
of $B$ on $B$ itself. 

It is proven in \cite{mon} that the $\s$-functor sending $E$ to $\VV(E)$ is fully faithful when
restricted to bounded above $k$-dg-modules $E$, and in particular for perfect
$k$-dg-modules $E$.

\section{Infinitesimal derived foliations}

As before, we fix $k$ a base simplicial commutative ring. \\

We let $\Ga$ be the additive group scheme over $k$, and consider
$G_0:=\Omega_0\Ga$ its derived loop scheme based at the origin. This is
a derived affine group scheme over $k$, and as an affine derived scheme
is given by $G_0\simeq \Spec\, k[\epsilon_{-1}]$, where $k[\epsilon_{-1}]$
is the free commutative 
simplicial $k$-algebra over one generator $\epsilon$ in homotopical degree $1$. The standard 
$\Gm$-action on $\Ga$ of weight $1$ induces a $\Gm$-action on $G_0$ corresponding to 
the grading on $k[\epsilon_{-1}]$ for which $\epsilon_{-1}$ is of weight $-1$. We then form the semi-direct
product $\cH_\pi:=\Gm \ltimes G_0$, which is a derived affine group scheme over $k$. \\

Before going further, let us mention that representations of the group $\cH_\pi$ are
precisely the graded mixed complexes for which the mixed structure is of cohomological degree $1$ (as 
opposed of the usual conventions taken for instance in \cite{ptvv,derfol,derfolbook} where the cohomological degree is $-1$). More precisely, 
we have an equivalence of symmetric monoidal $\s$-categories 
$$\QCoh(B\cH_\pi) \simeq +\medg_k,$$\
where $+\medg_k$ is the symmetric monoidal $\s$-category of graded mixed complexes with 
mixed structure of cohomological degree $1$. The $\s$-category of usual graded mixed complexes
with mixed structure of degree $-1$ is denoted by $\medg_k$, and is
itself given by representations of the group stack $\cH=\Gm\ltimes BKer$ of \cite{mrt}. 
It is well known that there
exists a symmetric monoidal equivalence of $\s$-categories
$$\mathsf{RS} : \medg_k \simeq +\medg_k$$
called the \emph{red-shift self equivalence}. The $\s$-functor $\mathsf{RS}$ sends a graded complex
$\oplus_n E^{(n)}$ to $\oplus_n E^{(n)}[-2n]$, and thus transforms a mixed structure $E^{(n)} \to E^{(n+1)}[-1]$ 
of cohomological degree $-1$ to a mixed structure $E^{(n)}[-2n] \to E^{(n+1)}[-2n-1]\simeq E^{(n+1)}[-2n+2][1]$
of cohomological degree $1$. To summarize, we have a commutative square of symmetric monoidal equivalences
$$\xymatrix{
\QCoh(B\cH) \ar[r]^-\sim \ar[d]_-\sim & \QCoh(B\cH_\pi) \ar[d]^-\sim \\
\medg_k \ar[r]_-{\mathsf{RS}} & +\medg_k.
}$$

It is proven in \cite{mrt} that $\QCoh(B\cH)$ is equivalent, as a symmetric monoidal $\s$-category, to 
the $\s$-category $\widehat{\fdg_k}$ of complete filtered complexes. This equivalence is 
given by the Tate realization $|-|^t$ sending a graded mixed complex $\oplus_n E^{(n)}$
to $\oplus_{i \leq 0}\prod_{p\geq i}E^{(p)}[-2p]$ endowed with the total 
differential, sum of the mixed structure and the cohomological differential. Similarly, 
$\QCoh(B\cH_\pi)$ is equivalent to $\widehat{\fdg_k}$ by means of the corresponding 
Tate realization $|-|^t$, sending $\oplus_n E^{(n)}$ to 
$\mathrm{colim}_{i \leq 0}(\prod_{p\geq i}E^{(p)})$ with the total differential. \\

For any derived stack $X \in \dSt_k$, we can form $\cL^{gr}_\pi(X/k):=\uMap(G_0,X)$, the derived
Hom stack (relative to $k$) from $G_0$ to $X$. It comes equipped with a natural action of $\cH$, where $G_0$ acts on itself
by translations and $\Gm$ acts by the grading on $G_0$. Together with this
$\cH_\pi$-action, $\cL^{gr}_\pi(X)$ is the cofree $\cH_\pi$-equivariant derived stack
generated by $X$ endowed with the trivial $\Gm$-action.

\begin{df}\label{dinfol}
Let $X=\Spec\, A$ be a derived affine scheme over $k$. A (perfect) \emph{infinitesimal derived foliation
on $X$ (relative to $k$)} consists of a $\cH_\pi$-equivariant derived stack $\F$ 
together with a $\cH_\pi$-equivariant morphism $p : \F \to X$, such that, 
as a $\Gm$-equivariant stack over $X$, $\F$ is a perfect  linear derived affine stack (i.e. of the
form $\VV(E)$ for $E$ a perfect complex on $X$).
\end{df}

The infinitesimal derived foliations on $X$ (relative to $k$) form an 
$\s$-category $\Fol^\pi(X/k)$, which by definition is a full sub-$\s$-category 
of $(\cH_\pi-\dSt_k)/\cL^{gr}_\pi(X/k)$, the $\cH_\pi$-equivariant derived stacks over
$\cL^{gr}_\pi(X/k)$.  The construction $X \mapsto \Fol^\pi(X/k)$ can be promoted to 
an $\s$-functor $\Fol^\pi(-/k) : \dAff_k^{op} \to \scat$, by
defining pull-backs as follows: for a morphism $f : X \to Y$ of affine derived $k$-schemes, 
we define $f^*$ to be induced by the base change
$$(\cH_\pi-\dSt_k)/\cL^{gr}_\pi(Y/k) \longrightarrow (\cH_\pi-\dSt_k)/\cL^{gr}_\pi(X/k),$$
along the induced morphism $\cL^{gr}_\pi(f/k): \cL^{gr}_\pi(X/k) \to \cL^{gr}_\pi(Y/k)$. As the
$\s$-functor $\Fol^\pi(-/k)$ clearly is a hyperstack for the fpqc topology, it 
can be uniquely extended to a colimit preserving $\s$-functor
$$\Fol^\pi(-/k) : \dSt_k \to \scat^{op}.$$

Let $X$ be a derived affine $k$-scheme, and $\F \in \Fol^\pi(X/k)$ be an infinitesimal derived foliation
on $X$. By \cite{mon}, the linear stack $\F \to X$ is of the form $\VV(\LL_{\F/k}[-1])$, for some uniquely 
defined perfect complex $\LL_{\F/k}$ on $X$. 

\begin{df}\label{dcot}
With the notations above, the complex $\LL_{\F/k}$ is called the \emph{cotangent complex
of $\F$ (relative to $k$)}. It is also simply denoted by $\LL_\F$ if $k$ is clear from the context.
\end{df}

When $X$ is a general derived stack and $\F \in \Fol^\pi(X/k)$, it is also 
possible to define the cotangent complex $\LL_{\F/k}$, as a quasi-coherent
complex over $X$. This construction is not straightforward, as
cotangent complexes in the sense above are not stable by pull-backs
and thus do not glue naively on stacks. 
When $X$ is moreover a derived Artin stack, the cotangent complex
$\LL_{\F/k}$ is a perfect complex over $X$. This global aspects of 
cotangent complexes will not be studied in details in this note, and we refer to \cite{derfolbook} for
more on the subject. However, the special case where $X$ is a derived Deligne-Mumford
stack presents no difficulties, as the cotangent complex of Definition \ref{dcot} are
stable by \'etale pull-backs, and thus will glue globally on $X$. This applies in particular
to the case where $X$ is any derived scheme. 

Let $X$ be a derived affine scheme of finite presentation over $k$.
The $\s$-category $\Fol^\pi(X/k)$ admits a final
object, namely $\F=\cL^{gr}_\pi(X/k)$. It is denoted by $*_{X/k}$, and its cotangent 
complex is naturally equivalent to $\LL_{X/k}$, the cotangent complex of $X$. Indeed, 
it is easy to see, simply by contemplating the functors of points and using the
very definition of the cotangent complex, that $\uMap(G_0,X) \simeq \VV(\LL_{X/k}[1])$.
Similarly, the initial object in $\Fol^\pi(X/k)$ exists, is denoted by $0_{X/k}$, and
its cotangent complex is $0$. It corresponds to the
constant-map morphism $\F=X \longrightarrow \cL^{gr}_\pi(X/k)$. 

\section{Infinitesimal cohomology}

Let $X$ be a derived Artin stack of finite presentation over $k$,  
and $\F \in \Fol^\pi(X/k)$ be an infinitesimal derived foliation on $X$. We define
the $\s$-category of quasi-coherent complexes along $\F$ as follows. 
We first consider $\QCoh([\F/\cH_\pi])$, the $\s$-category of quasi-coherent complexes
on the quotient stack $[\F/\cH_\pi]$ (which is also, by definition, the $\s$-category
of $\cH_\pi$-equivariant quasi-coherent complexes on $\F$). The $\s$-category
$\QCoh^\F(X)$ is defined as the full sub-$\s$-category of $\QCoh([\F/\cH_\pi])$
whose objects $E$ are graded free: the pull-back of $E$ on $[\F/\Gm]$ is of the form
$p^*(E_0)$ where $p : [\F/\Gm] \to X$ is the natural projection. 

\begin{df}\label{dinf}
With the notations above, 
the \emph{$\s$-category of quasi-coherent complexes along $\F$} is $\QCoh^\F(X)$. For $E \in 
\QCoh^\F(X)$, \emph{the infinitesimal cohomology of $X$ along $\F$ with coefficients in $E$} 
is defined to be
$$\Cinf^*(\F,E):=q_*(E) \in \dg_k$$
where $q : [\F/\cH_\pi] \longrightarrow \Spec\, k$ is the structural map.
\end{df}

The definition above will be justified by our next result, stating that 
when $\F=*_{X/k}$ and $E=\OO_X$, $\Cinf^*(\F,E)$ coincides with 
(completed derived) infinitesimal cohomology. For this, recall that 
if $X \to \Spec\, k$ is a smooth affine scheme, we have the infinitesimal
stack $X_{\mathrm{inf}}$ (relative to $k$), defined by the functor sending $A\in k-\sCR$ to 
$X(A_{red})$ (where $A_{red}=\pi_0(A)_{red}$ by definition). In characteristic
zero, the functor $X_{\mathrm{inf}}$ is often denoted by $X_{DR}$ (\cite{sim}). We prefer 
to avoid using the notation $X_{DR}$ in non-zero characteristic, and use $X_{\mathrm{inf}}$ instead, 
as the cohomology of $X_{\mathrm{inf}}$ computes infinitesimal cohomology and not de Rham cohomology. 
As any arbitrary derived stack, $X_{\mathrm{inf}}$ possesses a cohomology with coefficients in the structure sheaf $\OO$,
denoted by $C^*(X_{\mathrm{inf}},\OO)$.

\begin{thm}\label{compinf}
Let $X \to \Spec\, k$ be a smooth affine scheme (and assume $k$ is discrete
for simplicity). There exists an equivalence in $\dg_k$
$$\phi_X :  C^*(X_{\mathrm{inf}},\OO) \simeq  \Cinf^*(*_{X/k},\OO_X).$$
\end{thm}

\textit{Proof.} We first define the morphism $\phi_X$. For this we use the
canonical morphism $p : X \to X_{\mathrm{inf}}$, induced by $X(A) \to X(A_{red})$ coming from
the canonical projection $A \to A_{red}$. It induces a natural morphism of $\cH_\pi$-equivariant
derived stacks $\cL^{gr}_\pi(X) \to \cL^{gr}_\pi(X_{\mathrm{inf}})$. 
By the very definition of 
$X_{\mathrm{inf}}$, and because $G_0=\Omega_0\Ga$ is a derived affine scheme with trivial 
reduced sub-scheme, we see that the canonical morphism $X_{\mathrm{inf}} \to \cL^{gr}_\pi(X_{\mathrm{inf}})$
is in fact an ($\cH_\pi$-equivariant) equivalence of derived stacks. We thus obtain 
a natural $\cH_\pi$-equivariant morphism $\cL^{gr}_\pi(X) \to X_{\mathrm{inf}}$, 
or equivalently a morphism of derived stacks
$$q_X : [\cL^{gr}_\pi(X)/\cH_\pi] \to X_{\mathrm{inf}},$$
which is clearly functorial in $X$. It induces, by pull-back, a morphism between the cohomology
complexes with coefficients in $\OO$
$$\phi_X:=q_X^* :  C^*(X_{\mathrm{inf}},\OO) \to \Cinf(*_{X/k},\OO)$$
which is again functorial in $X$.

To prove that $\phi_X$ is an equivalence, we proceed by descent along $p : X \to X_{\mathrm{inf}}$.
As $X$ is smooth, the morphism $p$ is an epimorphism of derived stacks (infinitesimal
criterion for smoothness). The simplicial nerve $X_*$  of $p$, that comes
equipped with its augmentation $X_* \to X_{\mathrm{inf}}$, is an hypercovering. In particular, 
the canonical morphism
$$C^*(X_{\mathrm{inf}},\OO) \longrightarrow \lim_{[n] \in \Delta}\OO(X_n)$$
is a quasi-isomorphism. Note that $X_n$ is naturally isomorphic to the formal
completion of $X^n$ along the diagonal embedding $X \hookrightarrow X^n$, and 
thus $C^*(X_n, \OO)$ is the discrete algebra $\OO(X_n)$ of functions on the formal
scheme $X_n$.

On the other hand, the fact that $X$ is smooth also implies that 
$X \to [\cL^{gr}_\pi(X)/\cH_\pi]$ is an epimorphism, because $\cL^{gr}_\pi(X)\simeq \VV(\Omega_{X/k}^1[1])=B\T_X$
is the classifying stack of the tangent bundle of $X$.
We thus consider the commutative diagram of derived stacks with $\cH_\pi$-actions
$$\xymatrix{
X \ar[r]^-{id} \ar[d] & X \ar[d] \\
\cL^{gr}_\pi(X) \ar[r]_-{q_X} & X_{\mathrm{inf}}.
}$$
This diagram induces a morphism of simplicial objects on the nerves of the
vertical morphisms. The nerve of the morphism $X \to \cL^{gr}_\pi(X)$ can be described as the
simplicial object $[n] \mapsto \cL^{gr}_\pi(X/X^n)$, where $X$ sits inside $X^n$ via the diagonal embedding,
and $\cL^{gr}_\pi(X/X^n)$ is the relative mapping stack from $G_0$ to $X$ relative to $X^n$.
At the simplicial level $n$, the morphism induced on the nerve is
thus an $\cH_\pi$-equivariant morphism of stacks over $S$
$$q_{X,n} : \cL^{gr}_\pi(X/X^n) \longrightarrow X_n,$$
(where $X_n$ is the formal completion of $X^n$ along its diagonal).

\begin{lem}\label{laux1}
The morphism $q_{X,n}$ induces by pull-back a quasi-isomorphism of complexes
$$C^*(X_n,\OO) \simeq C^*([\cL^{gr}_\pi(X/X^n)/\cH_\pi],\OO).$$
\end{lem}

\textit{Proof of Lemma \ref{laux1}.} By descent we can easily restrict to the case where
$X=\Spec\, A$ is a smooth affine scheme over $S=\Spec\, k$. In this case,
$C^*(X_n,\OO)$ identifies canonically with $\hat{A}^{\otimes_k n}$, the formal
completion of $A^{\otimes_k n}$ along the  augmentation $A^{\otimes_k n} \to A$
induced by the multiplication map.
Similarly, $\cL^{gr}_\pi(X/X^n)$ is canonically identified with 
the smooth affine scheme $\Spec\, Sym_A((\Omega_{A/k}^1)^{n-1})$. Therefore, the morphism 
$\phi_{X,n}$ is a morphism of discrete commutative algebras
$$\phi_{A,n} : \hat{A}^{\otimes_k n} \longrightarrow |Sym_A((\Omega_{A/k}^1)^{n-1})|^t,$$
where the right hand side is the Tate realization of the graded mixed object $Sym_A((\Omega_{A/k}^1)^{n-1})$
induced by the $\cH_\pi$-action on $\cL^{gr}_\pi(X/X^n)$. The morphism $\phi_{A,n}$ is
compatible with the canonical augmentations to $A$, and thus is a morphism of complete
filtered commutative algebras. Therefore, in order to prove it is an isomorphism it is enough to show that 
it induces an isomorphism 
$$gr\phi_{A,n} : Gr^*(\hat{A}^{\otimes_k n})\longrightarrow Sym_A((\Omega_{A/k}^1)^{n-1})$$ on the associated graded objects.
Since $A$ is a smooth $k$-algebra, the left hand side is canonically isomorphic
to $Sym_A(I_n/I_n^2)$, where $I_n$ is the augmentation ideal. Moreover, 
we have the usual natural isomorphism of $A$-modules $I_n/I_n^2\simeq (\Omega_{A/k}^1)^{n-1}$, and thus
the morphism $gr\phi_{A,n}$ can be identified with a graded endomorphism
of the graded algebra $Sym_A((\Omega_{A/k}^1)^{n-1})$. As it is compatible with augmentations to $A$
it is thus determined by an $A$-linear endomorphism of $(\Omega_{A/k}^1)^{n-1}$. By functoriality in 
$A$ and in the simplicial direction $n$, we see moreover that it is determined by a functorial
endormorphism the $A$-module $\Omega_{A/k}^1$, and thus must be the multiplication by 
an element $\lambda \in k$. To conclude the lemma, we have to show that 
$\lambda=1$. 

By functoriality, it is easy to reduce to the case $A=k[T]$, a polynomial ring, and by base change we
can even assume that $k=\ZZ$. In order to check that $\lambda=1$, we can base change from $\ZZ$ to $\QQ$, and
thus even reduce to the case $A=\QQ[T]$ and $k=\QQ$. As we are now in characteristic zero, 
$C^*([\cL^{gr}_\pi(X/X^n)/\cH_\pi]$ identifies canonically with the 
completed de Rham complex $\CDR^*(X/X^n)$. Indeed, the red shift equivalence sends 
the graded mixed cdga $\OO(\cL^{gr}_\pi(X/X^n))$ to $\OO(\cL^{gr}(X/X^n))$. In such terms, 
the morphism $\phi_{X,n}$ is now the natural isomorphism relating 
derived de Rham cohomology of a closed embedding with functions on the formal 
completion (see \cite{cptvv}). \hfill $\Box$ \\

Lemma \ref{laux1} concludes the proof of Theorem \ref{compinf}, as we have a commutative diagram 
$$\xymatrix{
C^*(X_{\mathrm{inf}},\OO)  \ar[r]^-{\phi_X} \ar[d]_-{\sim} &  \Cinf(*_X,\OO_X)  \ar[d]^-{\sim} \\
 \lim_{[n]}\OO(X_n) \ar[r]_-{\sim} & \lim_{[n]}\Cinf(*_{X/X^n},\OO_X).}$$
\hfill $\Box$ \\

As a direct corollary of Theorem \ref{compinf}, we obtain that, for a general
derived Artin stack $X$, $\Cinf(*_{X/k},\OO)$ computes the \emph{completed derived 
infinitesimal cohomology of $X$} (relative to $k$). Indeed, it is easy to see the $\s$-functor
$X \mapsto \OO(\cL^{gr}_\pi(X/k))$, from affine derived schemes to 
graded mixed complexes, is equivalent to its extension by sifted colimits from 
smooth $X$'s. For $A \in \sCR$, where each $A_n$ is a smooth algebra over $k$, we thus have
$$\Cinf(*_X,\OO_X) \simeq \mathrm{colim}_{[n]} C^*((Spec\, A_n)_{\mathrm{inf}},\OO),$$
where the colimit is computed in the $\s$-category of \emph{complete} filtered complexes. The right
hand side of this equivalence is, by definition, the completed derived infinitesimal 
cohomology of $X$, and is denoted by $\mathbb{R}\hat{C}^*(X_{\mathrm{inf}},\OO)$. This notion is
extended to any derived Artin stack locally of finite presentation over $k$ by the usual gluing formula 
$$\mathbb{R}\hat{C}^*(X_{\mathrm{inf}},\OO):=\lim_{U \to X} \mathbb{R}\hat{C}^*(U_{\mathrm{inf}},\OO)$$
where the limit runs over affine $U$ locally of finite presentation.

\begin{cor}\label{c1}
For a general derived Artin stack $X$ of finite presentation over $k$, 
$\Cinf(*_{X/k},\OO)$ recovers the completed derived infinitesimal 
cohomology of $X$ relative to $k$
$$\Cinf(*_{X/k},\OO) \simeq \mathbb{R}\hat{C}^*(X_{\mathrm{inf}},\OO).$$
\end{cor}

\section{Integrability by formal groupoids}

We believe that infinitesimal derived foliations are formally integrable in a very strong sense
which we do not investigate in this paper in full generality. 
Rather, we prove here a special case, 
showing that the $\s$-category of smooth infinitesimal derived foliations over a $k$-scheme $X$
is in fact \emph{equivalent} to the category of \emph{formally smooth formal groupoids} over $X$. \\

Let $X$ be a smooth $k$-scheme and $\F\in \Fol^\pi(X/k)$. We say that $\F$ is \emph{smooth} if 
$\LL_\F$ is a vector bundle over $X$. In this case, we consider $X \times_\F X$, endowed
with its natural $\cH_\pi$-action and with its canonical projection to $X$. Note that, as a graded
derived stack, $X \times_\F X$ is of the form $\VV(\LL_{F})$, and thus is a the total 
space of the tangent bundle $\T_\F$ of $\F$. Also, the natural projection
$X\times_\F X \to X$ is $\cH_\pi$-equivariant. As a result, 
$X \times_\F X$ is a new object in $\Fol^\pi(X)$, that will be simply denoted by $\Omega_X\F$, and
called the \emph{loop space} of $\F$. Being the fiber product of the canonical morphism
$0_X \to \F$, the loop space $\Omega_X\F$ comes equipped with a natural structure of a groupoid object in $\Fol^\pi(X)$, acting
on the object $0_X$. 

\begin{prop}\label{p2}
The full sub-$\s$-category of $\Fol^\pi(X/k)$ formed by objects $\F$ whose cotangent complex
is of the form $V[1]$ for $V$ a vector bundle on $X$, is naturally equivalent to the 
category of formal schemes $Z$ with an identification $Z_{red} \simeq X$, 
which are locally equivalent on $X$ to $X\times \widehat{\mathbb{A}}^d$ where $d=rank(V)$.
\end{prop}

\textit{Proof.} To each $\F$ as in the proposition, we associate the sheaf 
of infinitesimal cohomology $\Cinf(\F)$ on $X$. By assumptions on the cotangent complex of $\F$, 
$\Cinf(\F)$ is a sheaf of discrete complete filtered commutative algebras, augmented 
to $\OO_X$ and 
its associated graded object is isomorphic to $\OO_X[[t_1,\dots,t_d]]$. Therefore, 
$Z_\F:=\Spf(\Cinf(\F))$ is a formal scheme with a canonical identification 
$Z_{red}\simeq X$ which is locally equivalent to $X\times \widehat{\mathbb{A}}^d$.

The fact that 
the construction $\F \mapsto Z_\F$ induces the equivalence of the proposition simply follows from the
fact that the Tate realization induces an equivalence between graded mixed complexes and
complete filtered complexes (see \cite[Prop. 1.3.1]{derfol2}). \hfill $\Box$ \\

By Proposition \ref{p1}, for any smooth infinitesimal derived foliation $\F$, its loop space 
$\Omega_X \F$ can be realized as a groupoid object in formal schemes acting on $X$
$$Z_{\Omega_X\F} \longrightarrow X\times X.$$

The construction $\F \mapsto Z_{\Omega_X\F}$ provides an functor from the
category of smooth infinitesimal derived foliations on $X$, and the category of 
formally smooth formal groupoids over $X$. We can produce a functor in the other direction, as follows.
Starting from a formally smooth formal groupoid $G \rightarrow X\times X$, we can consider its nerve
$G_* \to X$, which is a simplicial object in formal schemes. Taking functions on $G_*$
provides a cosimplicial complete filtered commutative algebra $\OO(G_*)$ whose associated graded
is the cosimplicial algebra of functions on the simplicial scheme $BV$, where 
$V=e^*(\Omega_{G/X}^1)$ is the vector bundle of invariant relative forms on $G$.
Using the equivalence between graded mixed complexes and complete filtered complexes of 
\cite[Prop. 1.3.1]{derfol2}, 
we find that $\OO(G_*)$ can be realized as a cosimplicial object in the category of
graded commutative algebras endowed with a compatible action on $G_0$. Passing then to $Spec$, we get 
a simplicial diagram of affine schemes endowed with a $\cH_\pi$-action, whose underlying 
simplicial diagram (obtained by forgetting the action) is the usual simplicial scheme $BV$. 
Taking geometric realization we get a $\cH_\pi$-action on the derived stack $\VV(V[1])$, 
which defines a smooth infinitesimal derived foliations $\F$ on $X$. 

We can then subsume the previous discussion in the following corollary.

\begin{cor}\label{cint}
The category of smooth infinitesimal derived foliations on $X$ is equivalent to the
category of formally smooth formal groupoids on $X$. 
\end{cor}

We note that the above corollary marks a major difference between 
derived foliations of \cite{derfol,derfolbook} and infinitesimal derived foliations. Indeed,
any finite dimensional Lie algebra $L$ over a field $k$, defines a 
smooth derived foliation over $k$, by considering $Sym(L^\vee[1])$ endowed with 
the mixed structure coming from the Chevalley-Eilenberg differential. However, 
it is not true that any Lie algebra can be integrated to a smooth formal
scheme in general. On the contrary, the above corollary, when applied to $X=\Spec\, k$ with $k$ a field, 
states that smooth infinitesimal derived foliations over $k$ form a category equivalent to 
smooth formal groups over $k$. Any such smooth infinitesimal derived foliation
is of the form $\VV(V[-1])$ for $V$ a finite dimensional $k$-vector space $V$. Therefore, the data
of an infinitesimal derived foliation structure on $\VV(V[-1])$ seems to be related to that of a \emph{partition
Lie algebra structure} in the sense of \cite{bracanu}. The precise comparison is out of the scope of this note, 
but some recent results of J. Fu (\cite{jfu}) make us strongly believe this is possible.

\section{Towards a comparison between infinitesimal derived foliations and derived foliations}

\subsection{Red shift}
We have already mentioned the red-shift equivalence $\mathsf{RS} : \medg_k \simeq +\medg_k$. We will now lift 
this to an endofunctor on the level of derived affine stacks endowed with actions of 
$\cH$ and $\cH_\pi$. For this, we use our cosimplicial-simplicial models for graded derived affine stacks.
We treat the absolute case of derived affine stacks over $\ZZ$, the realtive situation 
follows from considering graded affine stacks over a given affine derived scheme $X=\Spec\, k$ (trivially 
graded). For technical reasons we will actually invert $2$, i.e. most of the time work over $\ZZ[1/2]$.\\

We denote by $\ZZ[u]$ the cosimplicial commutative ring obtained by denormalization 
from the free commutative algebra over one generator in homotopical degree $0$ (endowed with the $0$
differential and considered as a non-negatively graded  commutative dg-algebra). In \cite{mrt},
it is proven that
$\Spec^\Delta(\ZZ[u])\simeq BS_{gr}^{1}$, the classifying stack of the \emph{graded circle}, 
also called the \emph{graded infinite projective space}. 
We consider $\ZZ[u]$ as an object in $csCR^{gr}$ constant in the simplicial direction where $u$ is of weight $1$.
Similarly, we denote by $\ZZ[v]$ the free simplicial commutative ring generated by 
one variable in homotopical degree $2$. We consider $\ZZ[v]$ as an object in $csCR^{gr}$
constant in the cosimplicial direction where $v$ is of weight $-1$.
Finally, we set 
$$\ZZ<u,v>:=\ZZ[v] \times_\ZZ \ZZ[u] \in csCR^{gr},$$
which is a cosimplicial-simplicial commutative $\ZZ$-graded ring.

For two commutative graded rings $A$ and $B$, we define their \emph{convolution product} $A\odot B$, which is 
a new commutative graded ring with 
$$(A\odot B)^{(n)}=A^{(n)}\otimes B^{(n)}$$
with the natural componentwise multiplication. This can be extended to a convolution product
for two objects in $csCR^{gr}$ by applying this construction levelwise, both in the simplicial and cosimplicial directions.

\begin{df}
The red shift endofunctor $\mathsf{RS} : csCR^{gr} \to csCR^{gr}$ is defined by
$$\mathsf{RS}(A):= A \odot (\ZZ<u,v>).$$
\end{df}

It is easy to see that $Tot^{\pi}(\mathsf{RS}(A))$ is a graded complex which is 
naturally quasi-isomorphic to $\mathsf{RS}(Tot^{\pi}(A))$. Therefore, $\mathsf{RS}$ preserves weak equivalence 
in $csCR^{gr}$, and thus induces a well defined $\s$-functor
$\mathsf{RS} : \csCR^{gr} \to \csCR^{gr}$ covering the red-shift self equivalence by the functor $Tot^\pi$
$$\xymatrix{
\csCR^{gr} \ar[r]^-{\mathsf{RS}} \ar[d]_-{Tot^{\pi}} & \csCR^{gr} \ar[d]^-{Tot^\pi} \\
\mdg_k \ar[r]_-{\mathsf{RS}} & \mdg_k.}$$
We note however that $\mathsf{RS}$ does not induce a self equivalence of $\csCR^{gr}$, and for instance
does not preserve free objects, as opposed for instance to the corresponding situation with  
$E_\infty$-algebras. In fact the behaviour of $\mathsf{RS}$ with respect to 
free objects is quite subtle, and the authors do not claim to fully understand the situation.

Note also that $\mathsf{RS}(A)$ is only well behaved when $A$ is only $\ZZ_\geq 0$-graded (or 
$\ZZ_{\leq 0}$-graded). For instance, $\mathsf{RS}(\ZZ[t,t^{-1}])$ is equivalent to 
$\ZZ<v,u>$, for which $u$ is the image of $t$ and $v$ the image of $t^{-1}$. However, by definition
$uv=0$, and thus the multiplicative structure on $\mathsf{RS}(\ZZ[t,t^{-1}])$ is somehow degenerate. 
Trying to use $\ZZ[u,u^{-1}]$ is not a solution here, as the presence of divided powers 
in the homotopy of commutative simplicial rings would force us to work over $\QQ$. \\

Using Proposition \ref{p2}, we can then consider $\mathsf{RS}$ as an endofunctor of the $\s$-category of graded derived affine stacks
$\mathsf{RS} : \dChAff^{gr} \to \dChAff^{gr}$. The convolution construction $-\odot B$
for $B\in csCR^{gr}$ is lax monoidal in the
sense that there is a natural morphism $(A \odot B) \otimes (A'\odot B) \to (A\otimes A')\odot B$. 
This lax monoidal structure is monoidal if and only if the multiplication maps
$B^{(p)} \otimes B^{(q)} \to B^{(p+q)}$ are completed quasi-isomorphisms of cosimplicial-simplicial
modules (i.e. becomes quasi-isomorphismes ater applying $Tot^\pi$). In the case
$B=\ZZ<u,v>$ this is never the case, as $uv=0$. However, if $A$ is $\ZZ_{\geq 0}$-graded, 
we have $A\odot \ZZ<u,v> \simeq A\odot \ZZ[u]$, and clearly $\ZZ[u]^{(p)} \otimes \ZZ[v]^{(q)} \to 
\ZZ[u]^{(p+q)}$ is a completed quasi-isomorphism as its image by $Tot^\pi$ is 
equivalent to the canonical isomorphism $\ZZ[-2p] \otimes \ZZ[-2q] \simeq \ZZ[-2p-2q]$.

\begin{lem}
Let $X$ and $Y$ be two graded derived affine stacks which are both $\ZZ_{\geq 0}$-graded. Then, 
the natural morphism $\mathsf{RS}(X\times Y) \to \mathsf{RS}(X)\times \mathsf{RS}(Y)$ is an equivalence. 
\end{lem}

Another special case where $\mathsf{RS}$ preserves products is given in the following Lemma. It follows from the easy observation that the image by $Tot^{\pi}$ of the multiplication
$\ZZ[v]^{(-1)} \otimes \ZZ[v]^{(-1)} \to \ZZ[v]^{(-2)}$  is equivalent to the multiplication by $2$
(because $\ZZ[v]$ has divided powers) 
$$\times 2 : \ZZ[2] \otimes \ZZ[2] \simeq \ZZ[4] \to \ZZ[4].$$

\begin{lem}\label{lem1/2}
Let $X$ and $Y$ be two graded derived affine stacks whose weights are concentrated in $[-1,0]$.
Then, the natural morphism $\mathsf{RS}(X\times Y) \to \mathsf{RS}(X)\times \mathsf{RS}(Y)$ is an equivalence
when restricted to $\Spec\, \ZZ[1/2]$. 
\end{lem}

Remind from \cite{mrt,derfol} 
the group object $BKer=S^{1}_{gr} \in \dChAff^{gr}$. As a graded derived affine stack 
its weights are concentrated in $[-1,0]$, so Lemma \ref{lem1/2} applies, and we deduce that 
that $\mathsf{RS}(BKer)$ is another group object in $\dChAff^{gr}$ (at least over $\Spec\, \ZZ[1/2]$). As a 
graded affine stack, it is clearly of the form $\Spec(\ZZ[1/2][\epsilon_{-1}])\simeq G_0$. 
As the group structure on $G_0$ is essentially unique (see \cite{mrt}), 
we moreover deduce that $\mathsf{RS}(BKer)\simeq G_0$ as group objects in $\dChAff^{gr}$.

We have thus seen that the image of the group $BKer$ by $\mathsf{RS}$ is the group $G_0$. We believe that 
the functor $\mathsf{RS}$ can also be promoted to a functor on the level of equivariant objects. This is not 
a formal statement, as $\mathsf{RS}$ does not commute with finite products in general. The precise construction
is outside of scope of this note and we leave this as an open question for the future.

\begin{quest}\label{q}
Can the $\s$-functor $\mathsf{RS}$ be extended to 
$$\mathsf{RS} : BKer-\dChAff^{gr+} \longrightarrow G_0-\dChAff^{gr+},$$
in such a way that the cofree object $X^{BKer}$ is sent to $\mathsf{RS}(X)^{G_0}$ (at least over $\Spec\, \ZZ[1/2]$) ?
\end{quest}

\begin{rmk}
\emph{
The above question is a coherence problem. It is indeed possible to prove easily that 
for any $X \in \dChAff^{gr+}$, there exists a canonical morphism 
$\mathsf{RS}(X^{BKer}) \to \mathsf{RS}(X)^{G_0}$ of graded derived affine stacks. Therefore, 
a $BKer$-action on $X$, induces a morphism $X \to X^{BKer}$, and therefore, via $\mathsf{RS}$, a morphism $\mathsf{RS}(X) \to \mathsf{RS}(X^{Bker}) \to \mathsf{RS}(X)^{G_0}$. The problem 
here is to control the higher coherences in order to promote the previous morphism
to an action of $G_0$ on $\mathsf{RS}(X)$.} 
\end{rmk}

\subsection{Infinitesimal structures on derived foliations}

We have now everything we need in order to compare infinitesimal derived foliations and 
our original definition of derived foliations,  assuming that we have a positive
answer to Question \ref{q}. 

Let $\F \to X$ be a derived foliation on a derived $k$-scheme $X$ in the sense of \cite{derfol}.
Remind that it consists of a linear derived stack $\VV(\LL_\F[1])$, for 
$\LL_\F$ a perfect complex on $X$, together with an action of the graded group $BKer$. 
By applying the red-shift functor, we get a graded derived affine stack $\mathsf{RS}(\VV(\LL_\F[1]))$
together with an action of $G_0$. The complex of functions on $\mathsf{RS}(\VV(\LL_\F[1]))$
is the graded complex $\mathsf{RS}(Sym_{\OO_X}(\LL_\F[1]))$ and thus receive a canonical 
morphism $\LL_\F[-1] \to \mathsf{RS}(Sym_{\OO_X}(\LL_\F[1]))$. By the universal property 
of the $Sym$ construction, we thus get a canonical morphism of graded stacks
$\mathsf{RS}(\VV(\LL_\F[1])) \to \VV(\LL_\F[-1]))$.

\begin{df}
Let $\F$ be a derived foliation on $X$ as above. An \emph{infinitesimal structure on $\F$}
consists of an extension of the $G_0$-action along the morphism $\mathsf{RS}(\VV(\LL_\F[1])) \to \VV(\LL_\F[-1]))$.
\end{df}

By definition, an infinitesimal structure on $\F$ determines an infinitesimal derived foliation
$\F^\pi$, given by the $G_0$-action on $\VV(\LL_\F[-1]))$. We can thus define the
infinitesimal cohomology of $\F$ by taking infinitesimal cohomology of $\F^\pi$ in the sense
of our Definition \ref{dinf}. 
Therefore, the first consequence of the data of an infinitesimal structure is the possibility to 
define infinitesimal cohomology. Note that with the definition above, there is a canonical 
morphism 
$$\Cinf(\F^\pi,\OO) \longrightarrow \CDR(\F,\OO)$$
from the infinitesimal cohomology of $\F^\pi$ to the derived de Rham cohomology of $\F$ introduced
in \cite{derfol}. This is a generalization of the well known comparison morphism between infinitesimal
and crystalline cohomology. This can be enhanced to a pull-back functor, from $\QCoh(\F^\pi)$
to $\QCoh(\F)$. Moreover, the differential operators construction of \cite{derfol2} can also be
performed for $\F$ and $\F^\pi$. This provides two different sheaves of filtered dg-algebras on $X$,
the crystalline differential operators $\D_\F$ along $\F$ and the Grothendieck differential operators
$\D^{\s}_{\F^\pi}$ along $\F^\pi$. These are two filtered dg-algebras, with a natural 
morphism $\D_\F \longrightarrow \D^{\s}_{\F^\pi}$, whose induced morphism on the associated graded objects is
$Sym(\T_\F) \longrightarrow Sym^{pd}(\T_\F)$, from the symmetric algebra (in the sense
of cosimplicial-simplicial models of our \S 1) to its PD-completion. 

We also believe that it should be possible to prove that a derived foliation that admits an
infinitesimal structure is automatically \emph{formally integrable}. In fact, it seems that 
when it exists, $\F^\pi$ completely recovers $\F$ as follows. The graded derived affine stack 
$\mathsf{RS}(\VV(\LL_\F[1]))$ can be obtained as the PD completion of $\VV(\LL_\F[-1]))$ along its
zero section, and it is very reasonable to expect that the $\mathsf{RS}$ construction
induces an equivalence between graded 
derived affine stacks and graded derived affine stacks with PD structures. In particular, the
$G_0=\mathsf{RS}(BKer)$-action on $\mathsf{RS}(\VV(\LL_\F[1]))$ should comes from a unique $BKer$-action
on $\VV(\LL_\F[1])$. What we are describing here is probably the existence of a
fogertful functor, from infinitesimal derived foliations to derived foliations, even though
we are not able, at the moment, to provide the details of this construction. This forgetful functor seems 
very closely related to the forgetful functor from partition Lie algebras to Lie algebras, and
recent results of J. Fu (\cite{jfu}) 
indicate that infinitesimal derived foliations could be essentially the same
thing as (perfect) partition Lie algebroids. The fact that partition Lie algebras are related to 
formal moduli problem then explains the fact that infinitesimal derived foliations have nice
formal integrability properties.

\bigskip
\bigskip
\bigskip
\bigskip
\bigskip

\noindent Bertrand To\"{e}n, {\sc IMT, CNRS, Universit\'e de Toulouse, Toulouse (France)}\\
Bertrand.Toen@math.univ-toulouse.fr \\

\medskip

\noindent Gabriele Vezzosi, {\sc DIMAI, Universit\`a di Firenze, Firenze (Italy)}\\
gabriele.vezzosi@unifi.it

\bigskip
\bigskip

\bibliographystyle{alpha}
\bibliography{Biblio.bib}

\newcommand{\etalchar}[1]{$^{#1}$}
\begin{thebibliography}{PTVV13}

\bibitem[BCN21]{bracanu}
Lukas Brantner, Ricardo Campos, and Joost Nuiten.
\newblock Pd operads and explicit partition lie algebras.
\newblock {\em Preprint arXiv:2104.03870}, 2021.

\bibitem[BM03]{MR2016697}
Clemens Berger and Ieke Moerdijk.
\newblock Axiomatic homotopy theory for operads.
\newblock {\em Comment. Math. Helv.}, 78(4):805--831, 2003.

\bibitem[CPT{\etalchar{+}}17]{cptvv}
Damien Calaque, Tony Pantev, Bertrand To\"{e}n, Michel Vaqui\'{e}, and Gabriele
  Vezzosi.
\newblock Shifted {P}oisson structures and deformation quantization.
\newblock {\em J. Topol.}, 10(2):483--584, 2017.

\bibitem[Eke87]{MR927978}
Torsten Ekedahl.
\newblock Foliations and inseparable morphisms.
\newblock In {\em Algebraic geometry, {B}owdoin, 1985 ({B}runswick, {M}aine,
  1985)}, volume~46 of {\em Proc. Sympos. Pure Math.}, pages 139--149. Amer.
  Math. Soc., Providence, RI, 1987.

\bibitem[Fu]{jfu}
Jiaqi Fu.
\newblock private communication.

\bibitem[Miy87]{MR927960}
Yoichi Miyaoka.
\newblock Deformations of a morphism along a foliation and applications.
\newblock In {\em Algebraic geometry, {B}owdoin, 1985 ({B}runswick, {M}aine,
  1985)}, volume~46 of {\em Proc. Sympos. Pure Math.}, pages 245--268. Amer.
  Math. Soc., Providence, RI, 1987.

\bibitem[Mon21]{mon}
Ludovic Monier.
\newblock A note on linear stacks, 2021.

\bibitem[MRT22]{mrt}
Tasos Moulinos, Marco Robalo, and Bertrand To\"{e}n.
\newblock A universal {H}ochschild-{K}ostant-{R}osenberg theorem.
\newblock {\em Geom. Topol.}, 26(2):777--874, 2022.

\bibitem[PTVV13]{ptvv}
Tony {Pantev}, Bertrand {To\"en}, Michel {Vaqui\'e}, and Gabriele {Vezzosi}.
\newblock {Shifted symplectic structures}.
\newblock {\em {Publ. Math., Inst. Hautes \'Etud. Sci.}}, 117:271--328, 2013.

\bibitem[Sim96]{sim}
Carlos Simpson.
\newblock Homotopy over the complex numbers and generalized de {R}ham
  cohomology.
\newblock In {\em Moduli of Vector Bundles (Taniguchi symposium December
  1994)}, volume 189 of {\em Lecture {N}otes in Pure and Applied Mathematics},
  pages 229--264. Dekker, 1996.

\bibitem[Toe06]{chaff}
Bertrand Toen.
\newblock Champs affines.
\newblock {\em Selecta Math. (N.S.)}, 12(1):39--135, 2006.

\bibitem[To{\"e}20]{derfol}
Bertrand To{\"e}n.
\newblock Classes caractéristiques des schémas feuilletés.
\newblock {\em Preprint arXiv:2008.10489}, 2020.

\bibitem[TV]{derfolbook}
Bertrand To\"{e}n and Gabriele Vezzosi.
\newblock Derived foliations.
\newblock {\em Book in preparation}.

\bibitem[TV08]{hagII}
Bertrand To\"{e}n and Gabriele Vezzosi.
\newblock Homotopical algebraic geometry. {II}. {G}eometric stacks and
  applications.
\newblock {\em Mem. Amer. Math. Soc.}, 193(902):x+224, 2008.

\bibitem[TV20]{derfol2}
Bertrand To\"{e}n and Gabriele Vezzosi.
\newblock Algebraic foliations and derived geometry: index theorems.
\newblock 2020.

\end{thebibliography}

\end{document}